\documentclass{qjmam}
\usepackage{amssymb}
\usepackage{amsmath}
\usepackage{graphicx}
\usepackage{amssymb}

\startpage{1} \yr{1999} \vol{52} \issue{1}

\newtheorem{definition}{Definition}[section]
\newtheorem{lemma}[definition]{Lemma}
\newtheorem{theorem}[definition]{Theorem}
\newtheorem{proposition}[definition]{Proposition}

\newtheorem{remark}[definition]{Remark}
\newtheorem{example}[definition]{Example}

\DeclareMathAlphabet\mathbit
    \encodingdefault\rmdefault\bfdefault\itdefault
\DeclareOldFontCommand{\bi}{\normalfont\bfseries\itshape}{\mathbit}

\newcommand{\be}{\begin{equation}}
\newcommand{\ee}{\end{equation}}

\def\fakebold#1{\relax\ifvmode\leavevmode\fi%
\ifmmode%
\setbox0=\hbox{$#1$}%
\else%
\setbox0=\hbox{#1}%
\fi%
\kern-.02em\copy0 \kern-\wd0%
\kern .04em\copy0 \kern-\wd0%
\kern-.0125em\raise.02em\box0%
}%

\renewcommand{\geq}{\geqslant}
\renewcommand{\leq}{\leqslant}

\begin{document}

\title[{about~global~minimizers~and~local~extrema~in phase~transition}]
{On D. Y. Gao and R. W. Ogden's paper ``Multiple solutions to
non-convex variational problems with implications for phase
transitions and numerical computation''}

\author[m.~d.~voisei \etal] {M. D. VOISEI }

\address{Towson University, Department of Mathematics, Towson, MD-21252, U.S.A.}

\extraauthor{C. Z\u{A}LINESCU\thanks{(zalinesc@uaic.ro)}}

\extraaddress{University {}``Al.I.Cuza'' Ia\c{s}i, Faculty of
Mathematics, 700506-Ia\c{s}i, Romania, and\\ Institute of Mathematics
Octav Mayer, Ia\c{s}i, Romania}

\received{\recd 19 February 2010. \revd 19 February 2010}

\maketitle

\eqnobysec

\begin{abstract}
In this note we prove that a recent result stated by D. Y. Gao and
R. W. Ogden on global minimizers and local extrema in a phase
transition problem is false. Our goal is achieved by providing a
thorough analysis of the context and result in question and
counter-examples.

\end{abstract}

\section{Introduction}

The optimization problem we have in focus is introduced on \cite[p.\
505]{Gao/Ogden:08} where one says {}``The primal variational problem
(1.1) for the soft device can be written in the form

$(\mathcal{P}_{s})$ : $\min\limits
_{u\in\mathcal{U}_{s}}{\displaystyle  \left\{
P_{s}(u)=\int_{0}^{1}\Bigl[\tfrac{1}{2}\mu
u_{x}^{2}+\tfrac{1}{2}\nu\left(\tfrac{1}{2}u_{x}^{2}-\alpha
u_{x}\right)^{2}\Bigr]dx-F(u)\right\} }$, $\quad$(3.2)''

\noindent where (see \cite[p.\ 501]{Gao/Ogden:08})

{}``$F(u)={\displaystyle \int_{0}^{1}}fudx+\sigma_{1}u(1)$
$\quad$(2.8)'',

\noindent and (see \cite[p.\ 505]{Gao/Ogden:08})

{}``$\mathcal{U}_{s}=\left\{ u\in\mathcal{L}(0,1)\mid
u_{x}\in\mathcal{L}^{4}(0,1), \ u(0)=0\right\} $. $\quad(3.1)$''

\smallskip
In Section 2 we explain the natural interpretation for the
definition of $\mathcal{U}_{s}$.

As mentioned on \cite[p.\ 498]{Gao/Ogden:08}, {}``$\mu$, $\nu$ and
$\alpha$ are positive material constants'', and {}``we focus mainly
on the case for which $\nu\alpha^{2}>2\mu$'' (see \cite[p.\
499]{Gao/Ogden:08}). Moreover (see \cite[p.\ 498]{Gao/Ogden:08}),
{}``To make the mixing of phases more dramatic, we introduce a
distributed axial loading (body force) $f\in\mathcal{C}[0,1]$ per
unit length of $I$''. These assumptions will be in force throughout
this article. Therefore, from

{}``$\sigma(x)={\displaystyle \int_{x}^{1}}f(s)ds+\sigma_{1}$
$\quad$(2.12)''

\noindent one obtains that $\sigma\in\mathcal{C}^{1}[0,1]$ and

{}``$F(u)={\displaystyle \int_{0}^{1}}\sigma(x)u_{x}dx$.
$\quad$(2.13)''

Furthermore (see \cite[p.\ 501]{Gao/Ogden:08}), one says {}``... we
obtain the Gao--Strang total complementary energy $\Xi(u,\zeta)$
(\textbf{16}) for this non-convex problem in the form

$\Xi(u,\zeta)=\cdots={\displaystyle \int_{0}^{1}\left[\tfrac{1}
{2}u_{x}^{2}(\zeta+\mu)-\alpha
u_{x}\zeta-\tfrac{1}{2}\nu^{-1}\zeta^{2}\right]dx-\int_{0}^{1}}fudx-\sigma_{1}u(1)$,
$\quad(2.7)$''.

In the text above (\textbf{16}) is our reference
\cite{Gao/Strang:89}.

In \cite[pp.\ 501, 502]{Gao/Ogden:08} one obtains {}``the so-called
\emph{pure complementary energy functional} (\textbf{7},
\textbf{17})

$P_{s}^{d}(\zeta)=-{\displaystyle
\frac{1}{2}\int_{0}^{1}\left(\frac{
(\sigma+\alpha\zeta)^{2}}{\mu+\zeta}+\nu^{-1}\zeta^{2}\right)}dx$,
$\quad$(2.14)

\noindent which is well defined on the dual feasible space

$\mathcal{S}_{a}=\left\{
\zeta\in\mathcal{L}^{2}\mid\zeta(x)+\mu\neq0,\
\zeta(x)\geq-\tfrac{1}{2}\nu\alpha^{2},\ \forall x\in[0,1]\right\}
.$''

\smallskip{}
 References (\textbf{7}, \textbf{17}) above are our references \cite{Gao:98}
and \cite{Gao:99}.

\medskip

Probably, by {}``well defined on ... $\mathcal{S}_{a}$'' the authors
of \cite{Gao/Ogden:08} mean that $P_{s}^{d}(\zeta)\in\mathbb{R}$ for
every $\zeta\in\mathcal{S}_{a}$. Note that \begin{equation}
\frac{(\sigma+\alpha\zeta)^{2}}{\mu+\zeta}+\nu^{-1}\zeta^{2}=\frac{\beta^{2}}
{\zeta+\mu}+2\alpha\beta+\alpha^{2}(\zeta+\mu)+\nu^{-1}\zeta^{2},
\label{pdsc}\end{equation}
 where (see \cite[p.\ 502]{Gao/Ogden:08})

{}``$\beta(x)=\sigma(x)-\mu\alpha,\quad\eta=(\nu\alpha^{2}-2\mu)^{3}/27\nu.$
$\quad(2.21)$''

\noindent and $\beta\in\mathcal{C}^{1}[0,1]$. Let us set $B_{0}:=
\{s\in[0,1]\mid\beta(s)=0\}$, $B_{0}^{c}:=[0,1]\setminus B_{0}$.

\medskip

Let $\zeta\in\mathcal{L}^{2}:=\mathcal{L}^{2}[0,1]$ and set
$E_{\zeta}:= \{x\in[0,1]\mid\zeta(x)+\mu=0\}$. In the sequel we use
the convention $0/0:=0$, which agrees with the convention
$0\cdot(\pm\infty):=0$ used in measure theory. With this convention
in mind, from (\ref{pdsc}), we obtain that
$P_{s}^{d}(\zeta)\in\mathbb{R}$ if and only if
$\frac{\beta^{2}}{\zeta+\mu}\in\mathcal{L}^{1}:=\mathcal{L}^{1}[0,1]$
which implicitly provides that $\frac{\beta^{2}}{\zeta+\mu}$ is
well-defined almost everywhere (a.e.\ for short), i.e.,
$E_{\zeta}\setminus B_{0}$ is negligible.

Consider \[ A_{1}:=\left\{
\zeta\in\mathcal{L}^{2}\mid\frac{\beta^{2}}{\zeta+
\mu}\in\mathcal{L}^{1}\right\} \subset A_{2}:=\left\{
\zeta\in\mathcal{L}^{2}\mid\zeta(x)+\mu\neq0\text{ for a.e.\ }x\in
B_{0}^{c}\right\} .\]
 The set $A_{1}$ is the greatest subset of $\mathcal{L}^{2}$ for
which $P_{s}^{d}(\zeta)\in\mathbb{R}$. Notice that
$\zeta\in\mathcal{L}^{2}$ makes $\frac{\beta^{2}}{\zeta+\mu}$ be
well-defined iff $\zeta\in A_{2}$. Also, note that
$\mathcal{S}_{a}\subset A_{2}$.

Denote by $\lambda$ the Lebesgue measure on $\mathbb{R}$. For
$\zeta\in A_{1}$ we have \begin{equation}
P_{s}^{d}(\zeta)=-\frac{1}{2}\int_{[0,1]\setminus E_{\zeta}}
\left(\frac{(\sigma+\alpha\zeta)^{2}}{\mu+\zeta}+\nu^{-1}\zeta^{2}
\right)dx-\frac{1}{2}\nu^{-1}\mu^{2}\lambda(E_{\zeta}).\label{pdsc2}\end{equation}

Notice that in the trivial case $\beta=0$ we have
$A_{1}=A_{2}=\mathcal{L}^{2}$ and so $P_{s}^{d}$ is well-defined on
$\mathcal{S}_{a}$ because in this case $P_{s}^{d}$ is well-defined
on $\mathcal{L}^{2}$.

\begin{proposition} \label{dom-gresit} If $\beta\neq0$ then $\mathcal{S}_{a}
\not\subset A_{1}$ and $P_{s}^{d}$ is not well-defined on
$\mathcal{S}_{a}$. \end{proposition}

\begin{proof} Because $\beta\neq0$ and $\beta\in\mathcal{C}^{1}[0,1]$,
there exist $\gamma>0$ and $0\leq a<b\leq1$ such that
$\beta^{2}(x)\geq\gamma$ for every $x\in[a,b]$. Consider
$\zeta(x):=x-a-\mu$ for $x\in(a,b)$ and $\zeta(x):=1-\mu$ for
$x\in[0,1]\setminus(a,b)$. Then $\zeta(x)>
-\mu>-\tfrac{1}{2}\nu\alpha^{2}$ for every $x\in[0,1]$ and
$\zeta\in\mathcal{L}^{2}$; hence $\zeta\in\mathcal{S}_{a}$. Note
that $\zeta\notin A_{1}$ since
$\frac{\beta^{2}}{\zeta+\mu}\ge\frac{\gamma}{\zeta+\mu}>0$ and
$\int_{0}^{1}\frac{1}{\zeta(x)+\mu}dx\ge\int_{a}^{b}\frac{dx}{x-a}=+\infty$.
\end{proof}


In the sequel $P_{s}^{d}$ is understood as being defined on $A_{1}$.

\medskip

Assume for the rest of this section that $\beta\neq0$. This yields
that $\lambda(B_{0}^{c})>0$ since $\beta$ is continuous.

\medskip

The next result is surely known. We give the proof for easy
reference.

\begin{lemma} Assume that $\sum_{n\geq1}\alpha_{n}<\infty$, where
$(\alpha_{n})_{n\geq1}\subset[0,\infty)$. Then there exists a
non-decreasing sequence $(\beta_{n})_{n\geq1}\subset(0,\infty)$ with
$\beta_{n}\rightarrow\infty$ and
$\sum_{n\geq1}\alpha_{n}\beta_{n}<\infty$. \end{lemma}

\begin{proof} Because the series $\sum_{n\geq1}\alpha_{n}$ is convergent,
the sequence $(R_{n})$ converges to $0$, where
$R_{n}:=\sum_{k=n+1}^{\infty}\alpha_{n}$. Hence there exists an
increasing sequence $(n_{k})_{k\geq1}\subset\mathbb{N}^{\ast}$ such
that $R_{n}<2^{-k}$ for all $k\geq1$ and $n\geq n_{k}$. Consider
$\beta_{n}:=1$ for $n\leq n_{1}$ and $\beta_{n}:=k$ for $n_{k}<n\leq
n_{k+1}$. Clearly, $(\beta_{n})$ is non-decreasing and
$\lim\beta_{n}=\infty$. Moreover, \begin{align*}
\sum_{p=1}^{n_{m+1}}\alpha_{p}\beta_{p} & =\sum_{p=1}^{n_{1}}
\alpha_{p}+\sum_{k=1}^{m}\sum_{p=n_{k}+1}^{n_{k+1}}\alpha_{p}\beta_{p}
\leq\sum_{p=1}^{n_{1}}\alpha_{p}+\sum_{k=1}^{m}k\sum_{p=n_{k}+1}^{n_{k+1}}\alpha_{p}\\
 & \leq\sum_{p=1}^{\infty}\alpha_{p}+\sum_{k=1}^{m}kR_{n_{k}}\leq
 \sum_{p=1}^{\infty}\alpha_{p}+\sum_{k=1}^{\infty}k2^{-k}<\infty.\end{align*}
 Therefore, the series $\sum_{n\geq1}\alpha_{n}\beta_{n}$ is convergent.
\end{proof}


Let us denote the algebraic interior (or core) of a set by
{}``$\operatorname*{core}$''.

\begin{proposition} Assume that $\beta\neq0$. Then $\operatorname*{core}A_{2}$
is empty. In particular,
$\operatorname*{core}A_{1}=\operatorname*{core}\mathcal{S}_{a}=\emptyset$.
\end{proposition}

\begin{proof} Let $\overline{\zeta}\in A_{2}$ be fixed. Then there
exists a sequence $(B_{n})_{n\geq1}$ of pairwise disjoint Lebesgue
measurable sets (even intervals) such that
$B_{0}^{c}=\cup_{n\geq1}B_{n}$ and $\lambda(B_{n})>0$ for $n\geq1$
(see e.g.\ \cite[p.\ 42]{roy-88}). We have that
$\sum_{n\geq1}\int_{B_{n}}\left\vert
\overline{\zeta}(x)+\mu\right\vert ^{2}dx=\int_{B_{0}^{c}}\left\vert
\overline{\zeta}(x)+\mu\right\vert ^{2}dx<\infty$, and so, from the
previous lemma, there exists a non-decreasing sequence
$(\beta_{n})_{n\geq1}\subset(0,\infty)$ with
$\beta_{n}\rightarrow\infty$ and \begin{equation}
\sum_{n\geq1}\beta_{n}\int_{B_{n}}\left\vert \overline{\zeta}(x)+
\mu\right\vert ^{2}dx<\infty.\label{star}\end{equation}
 Define $u:[0,1]\rightarrow\mathbb{R}$ by $u(x):=-\sqrt{\beta_{n}}(\overline{\zeta}(x)+\mu)$
for $x\in B_{n}$ and $u(x):=0$ for $x\in B_{0}$. From (\ref{star})
we have that $u\in\mathcal{L}^{2}$. Moreover, for every $\delta>0$
there exists a sufficiently large $N\geq1$ such that
$t=\beta_{N}^{-1/2}\in(0,\delta)$ and $\overline{\zeta}+tu\notin
A_{2}$; this happens because $B_{N}\subset \{x\in
B_{0}^{c}\mid\overline{\zeta}(x)+\beta_{N}^{-1/2}u(x)+\mu=0\}$ and
$\lambda(B_{N})>0$. We proved that
$\overline{\zeta}\not\in\operatorname*{core}A_{2}$. Hence
$\operatorname*{core}A_{2}=\emptyset$. \end{proof}


On page 502 of \cite{Gao/Ogden:08} it is said that {}``The
criticality condition with respect to $\zeta$ leads to the ... `dual
algebraic equation' (DAE) for ... (2.14) ..., namely

$\left(2\nu^{-1}\zeta+\alpha^{2}\right)(\mu+\zeta)^{2}=(\sigma-\mu\alpha)^{2}$.
$\quad$(2.16)''

\smallskip{}
 To our knowledge, one can speak about G\^{a}teaux differentiability of
a function $f:E\subset X\rightarrow Y$, with $X,Y$ topological
vector spaces, at $\overline{x}\in E$ only if $\overline{x}$ is in
the core of $E$. As we have seen above,
$P_{s}^{d}(\zeta)\in\mathbb{R}$ only for $\zeta\in A_{1}$ and
$\operatorname*{core}A_{1}=\emptyset.$

\medskip{}

\emph{So what is the precise critical point notion for} $P_{s}^{d}$
\emph{so that, when using that notion, one gets
\cite[(2.16)]{Gao/Ogden:08}, other than just formal computation?}

\medskip{}

Taking into account the comment (see \cite[p.\ 502]{Gao/Ogden:08})
{}``It should be pointed out that the integrand in each of
$P_{s}^{d}(\zeta)$ and $P_{h}^{d}(\zeta)$ has a singularity at
$\zeta=-\mu$, which explains the exclusion $\zeta\neq-\mu$ in the
definition of $\mathcal{S}_{a}$'', we must point out that there is
an important difference between the condition $\zeta\neq-\mu$ (as
measurable functions) and $\zeta(x)\neq-\mu$ a.e.\ on $[0,1]$ since
it is known that $\zeta\neq-\mu$ means that $\zeta(x)\neq-\mu$ on a
set of positive measure.

Alternatively, from the above considerations,
$\mathcal{L}^{2}\setminus\{-\mu\}$ is a (nonempty) open set, while
the set $A_{3}:=\left\{ \zeta\in\mathcal{L}^{2}
\mid\lambda(E_{\zeta})=0\right\} $ has, as previously seen, empty
core (in particular has empty interior).

The quoted text from \cite[p.\ 502]{Gao/Ogden:08} continues with:
{}``In fact, it turns out that, in general, $\zeta=-\mu$ does not
correspond to a critical point of either $P_{s}^{d}(\zeta)$ or
$P_{h}^{d}(\zeta)$. Exceptionally, we may have $\zeta(x)=-\mu$ for
some $x\in(0,1)$, but this is always associated with
$\sigma(x)=\mu\alpha$. It is therefore important to note that when
(2.16) holds, the integrand in (2.14) and (2.15) can be written as

$2\alpha(\sigma+\alpha\zeta)+\nu^{-1}\zeta(3\zeta+2\mu)$,
$\quad$(2.17)

\noindent and when $\zeta=-\mu$ (and $\sigma=\mu\alpha$) this
reduces to $\nu^{-1}\mu^{2}$, and the singularity in the integrand
is thus removed.''

This shows that the convention we used (namely $0/0=0$), our
interpretation for $P_{s}^{d}(\zeta)$, and formula (\ref{pdsc2}) are
in agreement with the authors of \cite{Gao/Ogden:08} point of view.

\section{Problem reformulation}

Every $u$ in $\mathcal{U}_{s}$ is represented by an absolutely
continuous function on $[0,1]$ with $u(0)=0$ and
$u_{x}\in\mathcal{L}^{4}(0,1)$. More accurately,
$\mathcal{U}_{s}=\left\{ u\in W^{1,4}(0,1)\mid u(0)=0\right\} $. In
a different notation, denoting by $\mathcal{L}^{p}$ the space
$\mathcal{L}^{p}[0,1]$, we have \[
u\in\mathcal{U}_{s}\Longleftrightarrow\exists v\in\mathcal{L}^{4},\
\forall x\in[0,1]:u(x)=\int_{0}^{x}v(t)dt.\]
 So, the problem $(\mathcal{P}_{s})$ above becomes \[
(\widehat{\mathcal{P}}_{s}):\quad\min_{v\in\mathcal{L}^{4}}\widehat{P}_{s}
 (v)=\int_{0}^{1}\left[\tfrac{1}{2}\mu
v^{2}+\tfrac{1}{2}\nu\left(\tfrac{1}{2}v^{2}-\alpha
v\right)^{2}-\sigma v\right]dx\]
 and $\Xi$ becomes \begin{equation}
\widehat{\Xi}(v,\zeta)=\int_{0}^{1}\left[\tfrac{1}{2}v^{2}(\zeta+\mu)-
\alpha v\zeta-\tfrac{1}{2}\nu^{-1}\zeta^{2}-\sigma
v\right]dx\quad(v\in\mathcal{L}^{4},\
\zeta\in\mathcal{L}^{2}).\label{def-Xi}\end{equation}

Note that $P_{s}(u)=\widehat{P}_{s}(v)$,
$\Xi(u,\zeta)=\widehat{\Xi}(v,\zeta)$, for
$u(x)=\int_{0}^{x}v(t)dt$, $x\in[0,1]$.

It is easy to see that $\widehat{P}_{s}$ and $\widehat{\Xi}$ are
Fr\'{e}chet differentiable and \begin{gather*}
d\widehat{P}_{s}(v)(h)=\int_{0}^{1}\left[\mu
v+\nu\left(\tfrac{1}{2}v^{2}-
\alpha v\right)(v-\alpha)-\sigma\right]hdx,\\
d\widehat{\Xi}(\cdot,\zeta)(v)(h)=\int_{0}^{1}\left[v(\zeta+\mu)-\alpha
\zeta-\sigma\right]hdx,\\
d\widehat{\Xi}(v,\cdot)(\zeta)(k)=\int_{0}^{1}\left[\tfrac{1}{2}v^{2}-
\alpha v-\nu^{-1}\zeta\right]kdx,\end{gather*}
 for $v,h\in\mathcal{L}^{4}$ and $\zeta,k\in\mathcal{L}^{2}$. Therefore,
\begin{gather}
\nabla \widehat{P}_{s}(v)=\mu v+\nu\left(\tfrac{1}{2}v^{2}-\alpha
v\right)
(v-\alpha)-\sigma\in\mathcal{L}^{4/3},\nonumber \\
\nabla
\widehat{\Xi}(\cdot,\zeta)(v)=v(\zeta+\mu)-\alpha\zeta-\sigma\in
\mathcal{L}^{4/3},\label{dzv}\\
\nabla \widehat{\Xi}(v,\cdot)(\zeta)=\tfrac{1}{2}v^{2}-\alpha
v-\nu^{-1}\zeta \in\mathcal{L}^{2}.\nonumber \end{gather}

Moreover, \begin{equation}
d^{2}\widehat{P}_{s}(v)(h,k)=\int_{0}^{1}\left[\mu+\nu\left(\tfrac{3}
{2}v^{2}-3\alpha
v+\alpha^{2}\right)\right]hkdx\quad(v,h,k\in\mathcal{L}^{4}).\label{d2ps}\end{equation}
 Hence $v\in\mathcal{L}^{4}$ is a critical point of $\widehat{P}_{s}$
if and only if \begin{equation} \mu
v+\nu\left(\tfrac{1}{2}v^{2}-\alpha v\right)(v-\alpha)-\sigma=0,
\label{cp-p}\end{equation}
 and $(v,\zeta)\in\mathcal{L}^{4}\times\mathcal{L}^{2}$ is a critical
point of $\widehat{\Xi}$ if and only if\begin{equation}
v(\zeta+\mu)-\alpha\zeta-\sigma=0,\quad\tfrac{1}{2}v^{2}-\alpha
v-\nu^{-1} \zeta=0.\label{cp-xi}\end{equation}

From the expression of $\widehat{\Xi}$ we observe that
$\widehat{\Xi}(v,\cdot)$ is concave on $\mathcal{L}^{2}$ for every
$v\in\mathcal{L}^{4}$; furthermore, $\widehat{\Xi}(\cdot,\zeta)$ is
convex (concave) for those $\zeta\in\mathcal{L}^{2}$ with
$\zeta\geq-\mu$ $(\zeta\leq-\mu)$.

\begin{lemma} \label{zv}Let $v\in\mathcal{L}^{4}$ and set \begin{equation}
\zeta_{v}:=\nu\left(\tfrac{1}{2}v^{2}-\alpha
v\right).\label{zetav}\end{equation}
 Then $\zeta_{v}\in\mathcal{L}^{2}$, $d\widehat{\Xi}(v,.)(\zeta)=0$
iff $\zeta=\zeta_{v}$, and \begin{equation}
\sup_{\zeta\in\mathcal{L}^{2}}\widehat{\Xi}(v,\zeta)=\widehat{\Xi}
(v,\zeta_{v})=\widehat{P}_{s}(v).\label{max-xi-v}\end{equation}

\end{lemma}

\begin{proof} The facts that for $(v,\zeta)\in\mathcal{L}^{4}\times\mathcal{L}^{2}$
we have $\zeta=\zeta_{v}$ iff $d\widehat{\Xi}(v,.)(\zeta)=0$ and
$\zeta_{v}\in\mathcal{L}^{2}$ are straightforward. Equality
(\ref{max-xi-v}) is due to the fact that every critical point
(namely $\zeta=\zeta_{v}$) of a concave function (namely
$\widehat{\Xi}(v,\cdot)$) is a global maximum point of that
function. \end{proof}


Consider the set \[ A_{0}:=\bigg\{
\zeta\in\mathcal{L}^{2}\mid\frac{\beta}{\zeta+\mu}\in
\mathcal{L}^{4}\bigg\} =\bigg\{
\zeta\in\mathcal{L}^{2}\mid\frac{\sigma-\alpha\mu}{\zeta+\mu}\in\mathcal{L}^{4}\bigg\}
,\]
 More precisely, $\zeta\in A_{0}$ iff $\zeta\in\mathcal{L}^{2}$,
$E_{\zeta}\subset B_{0}$, and
$\frac{\beta}{\zeta+\mu}\in\mathcal{L}^{4}([0,1] \setminus
E_{\zeta})$.

For $\zeta\in\mathcal{L}^{2}$ with $E_{\zeta}\subset B_{0}$ set
\begin{equation}
v_{\zeta}:=\frac{\sigma+\alpha\zeta}{\zeta+\mu}=\alpha+\frac{\beta}{\zeta+\mu}.
\label{vzeta}\end{equation}
 More precisely $v_{\zeta}(x)=\alpha+\frac{\beta(x)}{\zeta(x)+\mu}$
for $x\in[0,1]\setminus E_{\zeta}$ and $v_{\zeta}(x)=\alpha$ for
$x\in E_{\zeta}$. 
Notice that $\zeta\in A_{0}$ iff $v_{\zeta}\in\mathcal{L}^{4}$.

In the sequel $\chi_{E}$ denotes the characteristic function of
$E\subset[0,1]$, that is, $\chi_{E}(x)=1$ for $x\in E$ and
$\chi_{E}(x)=0$ for $x\in[0,1]\setminus E$.

\begin{lemma} \label{vz} For all $\zeta\in A_{0}$ and $v\in\mathcal{L}^{4}$
we have that
$d\widehat{\Xi}(\cdot,\zeta)(v_{\zeta}+\chi_{E_{\zeta}}v)=0$ and
$\widehat{\Xi}(v_{\zeta}+\chi_{E_{\zeta}}v,\zeta)=P_{s}^{d}(\zeta)$.

\end{lemma}

\begin{proof} According to (\ref{dzv}), we have \[
d\widehat{\Xi}(\cdot,\zeta)(v_{\zeta}+\chi_{E_{\zeta}}v)=(v_{\zeta}+
\chi_{E_{\zeta}}v)(\zeta+\mu)-\alpha\zeta-\sigma=\chi_{E_{\zeta}}
v(\zeta+\mu)=0\quad\forall\zeta\in A_{0},\ v\in\mathcal{L}^{4}.\]
Since $\zeta\in A_{0}$ we have $v_{\zeta}=\alpha$ and
$\sigma=\alpha\mu$ on $E_{\zeta}$. Taking into account
(\ref{def-Xi}), (\ref{pdsc2}) and using that outside $E_{\zeta}$ we
have $v_{\zeta}^{2}
(\zeta+\mu)=(\sigma+\alpha\zeta)v_{\zeta}=\frac{(\sigma+\alpha\zeta)^{2}}{\zeta+\mu}$,
we get\begin{align*}
\widehat{\Xi}(v_{\zeta}+\chi_{E_{\zeta}}v,\zeta)= &
-\tfrac{1}{2}\int_{[0,1]
\setminus E_{\zeta}}\left(\frac{(\sigma+\alpha\zeta)^{2}}{\mu+\zeta}+\nu^{-1}\zeta^{2}\right)dx\\
 & +\int_{E_{\zeta}}\left(\alpha\mu(\alpha+v)-\tfrac{1}{2}\nu^{-1}\mu^{2}-
 \sigma(\alpha+v)\right)dx=P_{s}^{d}(\zeta).\end{align*}
\end{proof}

In particular every $\zeta\in A_{0}$ is in the domain of
$P_{s}^{d}$, that is, $A_{0}\subset A_{1}$ (which can be observed
directly, too since $\beta\in\mathcal{L}^{\infty}$). The argument
above shows that $\widehat{\Xi}(\cdot,\zeta)$ has no critical points
if $\zeta\in\mathcal{L}^{2} \setminus A_{0}$ (due to the lack of
regularity) and $\widehat{\Xi}(\cdot,\zeta)$ has an infinity of
critical points of the form $v_{\zeta}+\chi_{E_{\zeta}}v$ with
$v\in\mathcal{L}^{4}$, if $\zeta\in A_{0}$ and
$\lambda(E_{\zeta})>0.$

Furthermore, for $\zeta\in A_{0}$, if $\zeta+\mu\geq0$
$(\zeta+\mu\leq0)$ then $v_{\zeta}$ is a global minimum (maximum)
point of $\widehat{\Xi}(\cdot,\zeta)$ because
$\widehat{\Xi}(\cdot,\zeta)$ is convex (concave) and $v_{\zeta}$ is
a critical point of $\widehat{\Xi}(\cdot,\zeta)$.
Hence\begin{equation} P_{s}^{d}(\zeta)=\left\{ \begin{array}{ccc}
\inf_{v\in\mathcal{L}^{4}}\widehat{\Xi}(v,\zeta) & \text{if} &
\zeta\in A_{0}
\text{ and }\zeta\geq-\mu,\\
\sup_{v\in\mathcal{L}^{4}}\widehat{\Xi}(v,\zeta) & \text{if} &
\zeta\in A_{0}\text{ and
}\zeta\leq-\mu.\end{array}\right.\label{pds}\end{equation}

\begin{theorem} \label{analysis} $\quad$\\

\vspace{-0.4cm}

{(i)} Let
$(\overline{v},\overline{\zeta})\in\mathcal{L}^{4}\times\mathcal{L}^{2}$
be a critical point of $\widehat{\Xi}$. Then
$\zeta_{\overline{v}}=\overline{\zeta}$,
$v_{\overline{\zeta}}=(1-\chi_{E_{\overline{\zeta}}})\overline{v}+
\alpha\chi_{E_{\overline{\zeta}}}\in\mathcal{L}^{4}$, $\overline{v}$
is a critical point of $\widehat{P}_{s}$, $\overline{\zeta}\in
A_{0}$,
$\widehat{P}_{s}(\overline{v})=\widehat{\Xi}(\overline{v},\overline{\zeta})=
P_{s}^{d}(\overline{\zeta})$,
$\left(2\nu^{-1}\overline{\zeta}+\alpha^{2}\right)(\mu+\overline{\zeta})^{2}=
(\sigma-\mu\alpha)^{2}$ (i.e. $\overline{\zeta}$ satisfies
\cite[(2.16)]{Gao/Ogden:08}), and
\begin{equation}
d^{2}\widehat{P}_{s}(\overline{v})(h,k)=3\int_{0}^{1}\left(\overline{\zeta}-
\rho\right)hkdx\label{d2ps-cr}\end{equation}
 for $h,k\in\mathcal{L}^{4},$ where \begin{equation}
\rho:=-\tfrac{1}{3}\left(\mu+\nu\alpha^{2}\right).\label{rho}\end{equation}
 If, in addition, $\overline{\zeta}\geq-\mu$ then \begin{equation}
\sup_{\zeta\in\mathcal{L}^{2}}\inf_{v\in\mathcal{L}^{4}}\widehat{\Xi}(v,\zeta)
=\inf_{v\in\mathcal{L}^{4}}\widehat{\Xi}(v,\overline{\zeta})=\widehat{\Xi}
(\overline{v},\overline{\zeta})=\widehat{P}_{s}(\overline{v})=\inf_{v\in
\mathcal{L}^{4}}\widehat{P}_{s}(v)=P_{s}^{d}(\overline{\zeta})=\sup_{\zeta\in
A_{0},\zeta\geq-\mu}P_{s}^{d}(\zeta).\label{equa}\end{equation}
 In particular $\overline{v}$ is a global minimum of $\widehat{P}_{s}$
on $\mathcal{L}^{4}$.

\medskip{}

{(ii)} If $v\in\mathcal{L}^{4}$ is a critical point of
$\widehat{P}_{s}$ then
$(v,\zeta_{v})\in\mathcal{L}^{4}\times\mathcal{L}^{2}$ is a critical
point of $\widehat{\Xi}$.

\medskip{}

{(iii)} Assume that $\zeta$ is a measurable solution of
$\left(2\nu^{-1}
\zeta+\alpha^{2}\right)(\mu+\zeta)^{2}=(\sigma-\mu\alpha)^{2}$ and
$v\in\mathcal{L}^{4}$. Then:

{(a)} $\zeta\in A_{0}$ and
$(v_{\zeta},\zeta)\in\mathcal{L}^{\infty}\times\mathcal{L}^{\infty}\subset
\mathcal{L}^{4}\times\mathcal{L}^{2}$. Moreover, \begin{equation}
\widehat{P}_{s}(v_{\zeta}+v\chi_{E_{\zeta}})=P_{s}^{d}(\zeta)+\tfrac{1}{8}
\nu\int_{E_{\zeta}}(v^{2}-\alpha^{2}+2\nu^{-1}\mu)^{2}dx\label{pps}\end{equation}
 and $(v_{\zeta}+v\chi_{E_{\zeta}},\zeta)$ is a critical point of
$\widehat{\Xi}$ iff
$\widehat{P}_{s}(v_{\zeta}+v\chi_{E_{\zeta}})=P_{s}^{d}(\zeta)$ iff
\begin{equation} v^{2}-\alpha^{2}+2\nu^{-1}\mu=0~\text{a.e.\ in
}E_{\zeta}.\label{valfa}\end{equation}
 In particular, $(v_{\zeta},\zeta)$ is a critical point of $\widehat{\Xi}$
iff $\lambda(E_{\zeta})=0$. 

{(b)} $v_{\zeta}+v\chi_{E_{\zeta}}$ is a critical point of
$\widehat{P}_{s}$ iff \[
v\left(v^{2}-\alpha^{2}+2\nu^{-1}\mu\right)=0~\text{a.e.\ in
}E_{\zeta}.\]

\end{theorem}

\begin{proof} (i) Assume that $(\overline{v},\overline{\zeta})\in
\mathcal{L}^{4}\times\mathcal{L}^{2}$ is a critical point of
$\widehat{\Xi}$. From (\ref{cp-xi}) we see that
$\overline{\zeta}=\zeta_{\overline{v}}$,
$v_{\overline{\zeta}}=(1-\chi_{E_{\zeta}})\overline{v}+\alpha
\chi_{E_{\overline{\zeta}}}\in\mathcal{L}^{4}$ which provides
$\overline{\zeta}\in A_{0}$, $\overline{v}$ is a critical point of
$\widehat{P}_{s}$, and
$\left(2\nu^{-1}\overline{\zeta}+\alpha^{2}\right)(\mu+\overline{\zeta})^{2}
=(\sigma-\mu\alpha)^{2}$. Note that
$v_{\overline{\zeta}}+\chi_{E_{\overline{\zeta}}}(\overline{v}-\alpha)=
\overline{v}$. The equality
$\widehat{P}_{s}(\overline{v})=\widehat{\Xi}(\overline{v},
\overline{\zeta})=P_{s}^{d}(\overline{\zeta})$ is a consequence of
Lemmas \ref{zv}, \ref{vz}.

Taking into account (\ref{d2ps}) and the second equation in
(\ref{cp-xi}) we obtain that for $h,k\in\mathcal{L}^{4},$ \[
d^{2}\widehat{P}_{s}(\overline{v})(h,k)=\int_{0}^{1}\left[\mu+\nu\left(3\nu^{-1}
\overline{\zeta}+\alpha^{2}\right)\right]hkdx=3\int_{0}^{1}
\left(\overline{\zeta}-\rho\right)hkdx.\]

Assume, in addition, that $\overline{\zeta}\geq-\mu$. Therefore
$\widehat{\Xi}(\cdot,\overline{\zeta})$ is convex and
$P_{s}^{d}(\overline{\zeta})=\inf_{v\in\mathcal{L}^{4}}\widehat{\Xi}(v,
\overline{\zeta})$ (see (\ref{pds})). Since $\overline{v}$ is a
critical point it yields that $\overline{v}$ is a global minimum
point of $\widehat{\Xi}(\cdot,\overline{\zeta})$. Similarly,
$\overline{\zeta}$ is a global maximum point for the concave
function $\widehat{\Xi}(\overline{v},\cdot)$. We get \[
\widehat{\Xi}(v,\overline{\zeta})\geq\widehat{\Xi}(\overline{v},\overline{\zeta})
\geq\widehat{\Xi}(\overline{v},\zeta)\quad\forall
v\in\mathcal{L}^{4},\ \forall\zeta\in\mathcal{L}^{2}.\]
 This implies that \[
\sup_{\zeta\in\mathcal{L}^{2}}\inf_{v\in\mathcal{L}^{4}}\widehat{\Xi}(v,\zeta)
\geq\inf_{v\in\mathcal{L}^{4}}\widehat{\Xi}(v,\overline{\zeta})=\widehat{\Xi}
(\overline{v},\overline{\zeta})=\sup_{\zeta\in\mathcal{L}^{2}}\widehat{\Xi}
(\overline{v},\zeta)\geq\inf_{v\in\mathcal{L}^{4}}\sup_{\zeta\in\mathcal{L}^{2}}
\widehat{\Xi}(v,\zeta).\]
 Since $\sup_{\zeta\in\mathcal{L}^{2}}\inf_{v\in\mathcal{L}^{4}}
 \widehat{\Xi}(v,\zeta)\leq\inf_{v\in\mathcal{L}^{4}}\sup_{\zeta\in\mathcal{L}^{2}}
 \widehat{\Xi}(v,\zeta)$
(this happens for every function $\widehat{\Xi}$), we obtain
together with (\ref{max-xi-v}) that\begin{equation}
\sup_{\zeta\in\mathcal{L}^{2}}\inf_{v\in\mathcal{L}^{4}}\widehat{\Xi}(v,\zeta)=
\inf_{v\in\mathcal{L}^{4}}\widehat{\Xi}(v,\overline{\zeta})=\widehat{\Xi}
(\overline{v},\overline{\zeta})=\widehat{P}_{s}(\overline{v})=\inf_{v\in
\mathcal{L}^{4}}\widehat{P}_{s}(v)=P_{s}^{d}(\overline{\zeta}).\label{minmax}\end{equation}
 In particular, $\overline{v}$ is a global minimum of $\widehat{P}_{s}$
on $\mathcal{L}^{4}$.

From (\ref{pds}) and (\ref{minmax}) we have
$$P_{s}^{d}(\overline{\zeta})=
\inf_{v\in\mathcal{L}^{4}}\widehat{\Xi}(v,\overline{\zeta})\le\sup_{\zeta\in
A_{0},\zeta\geq-\mu}\inf_{v\in\mathcal{L}^{4}}\widehat{\Xi}(v,\zeta)\le
\sup_{\zeta\in\mathcal{L}^{2}}\inf_{v\in\mathcal{L}^{4}}\widehat{\Xi}
(v,\zeta)=P_{s}^{d}(\overline{\zeta}).$$

The assertion (ii) follows directly from (\ref{cp-p}) and
(\ref{cp-xi}).

(iii) For given $\beta\in\mathcal{C}^{1}[0,1]$ relation
$\left(2\nu^{-1}\zeta+
\alpha^{2}\right)(\mu+\zeta)^{2}=\beta^{2}(x)$
$(=(\sigma(x)-\mu\alpha)^{2})$ is a polynomial equation in $\zeta$.

Let $\zeta:[0,1]\to\mathbb{R}$ be such that $\zeta(x)$ is a solution
of the previous equation for every $x\in[0,1]$, that is, $\zeta$ is
a solution of \cite[(2.16)]{Gao/Ogden:08}. Because $\beta^{2}$ is
bounded (being continuous) we have that $\zeta$ is bounded. If, in
addition, $\zeta$ is measurable then
$\zeta\in\mathcal{L}^{\infty}\subset\mathcal{L}^{2}$.

(a) Note that, due to \cite[(2.16)]{Gao/Ogden:08}, $E_{\zeta}\subset
B_{0}$ and $v_{\zeta}=\alpha+\beta/(\mu+\zeta)$ outside $E_{\zeta}$
whence
$(v_{\zeta}-\alpha)^{2}=2\nu^{-1}\zeta+\alpha^{2}\in\mathcal{L}^{\infty}([0,1]
\setminus E_{\zeta})$. Therefore
$v_{\zeta}\in\mathcal{L}^{\infty}\subset\mathcal{L}^{4}.$ This shows
that $\zeta\in A_{0}$.

Let $v\in\mathcal{L}^{4}$. Recall that
$v_{\zeta}+v\chi_{E_{\zeta}}=\alpha+v$, $\sigma=\alpha\mu$,
$\zeta=-\mu$ inside $E_{\zeta}$ and
$v_{\zeta}+v\chi_{E_{\zeta}}=v_{\zeta}$ outside $E_{\zeta}$, and so
\begin{align} \widehat{P}_{s}(v_{\zeta}+v\chi_{E_{\zeta}})= &
\int_{[0,1]\setminus E_{\zeta}} \left[\tfrac{1}{2}\mu
v_{\zeta}^{2}+\tfrac{1}{2}\nu\left(\tfrac{1}{2}v_{\zeta}^{2}
-\alpha v_{\zeta}\right)^{2}-\sigma v_{\zeta}\right]dx\notag\\
 & +\int_{E_{\zeta}}\left[\tfrac{1}{2}\mu(\alpha+v)^{2}+\tfrac{1}{2}\nu
 \left(\tfrac{1}{2}(\alpha+v)^{2}-\alpha(\alpha+v)\right)^{2}-
 \alpha\mu(\alpha+v)\right]dx.\label{r-phs}\end{align}

Taking into account that $\zeta(x)$ is a solution of the equation
\cite[(2.16)]{Gao/Ogden:08} and that for $x\in[0,1]\setminus
E_{\zeta}$ one has $\zeta(x)+\mu\neq0$, one gets \[ \tfrac{1}{2}\mu
v_{\zeta}^{2}+\tfrac{1}{2}\nu\left(\tfrac{1}{2}v_{\zeta}^{2}- \alpha
v_{\zeta}\right)^{2}-\sigma
v_{\zeta}=-\tfrac{1}{2}\frac{(\sigma+\alpha\zeta)^{2}}{\zeta+\mu}-\tfrac{1}{2}
\nu^{-1}\zeta^{2}\quad\text{on }[0,1]\setminus E_{\zeta}.\]
 A simple verification shows that \[
\tfrac{1}{2}\mu(\alpha+v)^{2}+\tfrac{1}{2}\nu\left(\tfrac{1}{2}(\alpha+v)^{2}-
\alpha(\alpha+v)\right)^{2}-\alpha\mu(\alpha+v)=\tfrac{1}{8}\nu(v^{2}-\alpha^{2}+
2\nu^{-1}\mu)^{2}-\tfrac{1}{2}\nu^{-1}\mu^{2}.\]
 Using the preceding equalities, from (\ref{r-phs}) and (\ref{pdsc2})
we obtain that (\ref{pps}) holds.

A direct computation shows that
$(v_{\zeta}+v\chi_{E_{\zeta}},\zeta)$ is a critical point of
$\widehat{\Xi}$ if and only if $v^{2}-\alpha^{2}+2\nu^{-1}\mu=0$
a.e.\ in $E_{\zeta}$. Therefore the mentioned equivalencies are
true. Moreover, because $\nu\alpha^{2}>2\mu$ the last equivalence
holds, too.

(b) Similarly, $v_{\zeta}+v\chi_{E_{\zeta}}$ is a critical point of
$\widehat{P}_{s}$ if and only if
$v(v^{2}-\alpha^{2}+2\nu^{-1}\mu)=0$ a.e.\ in $E_{\zeta}$.
\end{proof}

\medskip

Note the following direct consequences of the previous theorem:
\begin{itemize}
\item if $v\in\mathcal{L}^{4}$ is a critical point of $\widehat{P}_{s}$,
then $(v,\zeta_{v})$ is a critical point of $\widehat{\Xi}$,
$\zeta_{v}\in\mathcal{L}^{2}$ is a solution of
\cite[(2.16)]{Gao/Ogden:08}, and $\widehat{P}_{s}(v)=
\widehat{\Xi}(v,\zeta_{v})=P_{s}^{d}(\zeta_{v})$;
\item if $\zeta$ is a measurable solution of \cite[(2.16)]{Gao/Ogden:08}
and $v\in\mathcal{L}^{4}$ satisfies (\ref{valfa}) then
$\zeta=\zeta_{(v_{\zeta}+ v\chi_{E_{\zeta}})}$ and
$v_{\zeta}+v\chi_{E_{\zeta}}$ is a global minimum of
$\widehat{P}_{s}$ on $\mathcal{L}^{4}$;
\item it is possible $v_{\zeta}+v\chi_{E_{\zeta}}$ to be a critical point
of $\widehat{P}_{s}$ without $(v_{\zeta}+v\chi_{E_{\zeta}},\zeta)$
being a critical point of $\widehat{\Xi}$; such a situation happens
when $v=0$ and $\lambda(E_{\zeta})>0.$
\end{itemize}

\section{Discussion of \cite[Th.\ 3]{Gao/Ogden:08}}

Based on the above considerations we discuss the result in
\cite[Th.\ 3]{Gao/Ogden:08}; for completeness we also quote its
proof. Recall that

{}``$\beta(x)=\sigma(x)-\alpha\mu,\quad\eta=(\nu\alpha^{2}-2\mu)^{3}/27\nu.$
$\quad(2.21)$''

\medskip

{}``\textsc{Theorem 3}. (Global minimizer and local extrema) Suppose
that the body force $f(x)$ and dead load $\sigma_{1}$ are given and
that $\sigma(x)$ is defined by (2.12). Then, if $\beta^{2}(x)>\eta$,
$\forall x\in(0,1)$, the DAE (2.16) has a unique solution
$\overline{\zeta}(x)>-\mu$, which is a global maximizer of
$P_{s}^{d}$ over $\mathcal{S}_{a}$, and the corresponding solution
$\overline{u}(x)$ is a global minimizer of $P_{s}(u)$ over
$\mathcal{U}_{s}$,

$P_{s}(\overline{u})=\min\limits _{u\in\mathcal{U}_{s}}P_{s}(u)=\max
\limits
_{\zeta\in\mathcal{S}_{a}}P_{s}^{d}(\zeta)=P_{s}^{d}(\overline{\zeta}).\quad(3.9)$

If $\beta^{2}(x)\leq\eta$, $\forall x\in(0,1)$, then (2.16) has
three real roots ordered as in (3.5). Moreover,
$\overline{\zeta}_{1}(x)$ is a global maximizer of
$P_{s}^{d}(\zeta)$ over the domain $\zeta>-\mu$, the corresponding
solution $\overline{u}_{1}(x)$ is a global minimizer of $P_{s}(u)$
over $\mathcal{U}_{s}$ and

$P_{s}(\overline{u}_{1})=\min\limits
_{u\in\mathcal{U}_{s}}P_{s}(u)=\max \limits
_{\zeta>-\mu}P_{s}^{d}(\zeta)=P_{s}^{d}(\overline{\zeta}_{1}).\quad(3.10)$

For $\overline{\zeta}_{2}(x)$ and $\overline{\zeta}_{3}(x)$, the
corresponding solutions $\overline{u}_{2}(x)$ and
$\overline{u}_{3}(x)$ are, respectively, a local minimizer and a
local maximizer of $P_{s}(u)$,

$P_{s}(\overline{u}_{2})=\min\limits _{u\in\mathcal{U}_{2}}P_{s}(u)=
\min\limits
_{\overline{\zeta}_{3}<\zeta<-\mu}P_{s}^{d}(\zeta)=P_{s}^{d}
(\overline{\zeta}_{2})\quad(3.11)$

\noindent and

$P_{s}(\overline{u}_{3})=\max\limits _{u\in\mathcal{U}_{3}}P_{s}(u)=
\max\limits
_{-\tfrac{1}{2}\nu\alpha^{2}<\zeta<\overline{\zeta}_{2}}P_{s}^{d}(\zeta)
=P_{s}^{d}(\overline{\zeta}_{3}),\quad(3.12)$

\noindent where $\mathcal{U}_{j}$ is a neighborhood of
$\overline{u}_{j}$, for $j=2,3$.

\medskip{}

\emph{Proof.} This theorem is a particular application of the
general analytic solution obtained in (\textbf{7}, \textbf{14})
following triality theory.''

\bigskip{}

Note that (\textbf{7}, \textbf{14}) are our references \cite{Gao:98}
and \cite{Gao:00}.

\smallskip{}
 Before discussing the previous result let us clarify the meaning
of $\overline{\zeta}_{i}$ and $\overline{u}_{i}$ (as well as
$\overline{\zeta}$ and $\overline{u}$) appearing in the statement
above. Actually these functions are introduced in the statement of
\cite[Th.~2]{Gao/Ogden:08}:

\medskip{}

{}``\textsc{Theorem} 2. (Closed-form solutions) For a given body
force $f(x)$ and dead load $\sigma_{1}$ such that $\sigma(x)$ is
defined by (2.12), the DAE (2.16) has at most three real roots
$\overline{\zeta}_{i}(x)$, $i=1,2,3$, given by (2.22)--(2.24) and
ordered as

\medskip{}

$\overline{\zeta}_{1}(x)\geq-\mu\geq\overline{\zeta}_{2}(x)\geq
\overline{\zeta}_{3}(x)\geq-\tfrac{1}{2}\nu\alpha^{2}.\quad(3.5)$

\medskip{}

\noindent For $i=1$, the function defined by

\medskip{}

$\overline{u}_{i}(x)={\displaystyle
\int_{0}^{x}\frac{\sigma(s)+\alpha
\overline{\zeta}_{i}(s)}{\overline{\zeta}_{i}(s)+\mu}ds}\quad$ (3.6)

\medskip{}

\noindent is a solution of (BVP1). For each of $i=2,3$, (3.6) is
also a solution of (BVP1) provided $\overline{\zeta}_{i}$ is
replaced by $\overline{\zeta}_{1}$ for values of $s\in\lbrack0,x)$
for which $\overline{\zeta}_{i}(s)$ is complex.

For a given $t$ such that $\sigma_{1}$ is determined by $(3.3)_{3}$,
one of $\overline{u}_{i}(x)$, $i=1,2,3$, satisfies $(3.4)_{3}$ and
hence solves (BVP2). Furthermore,

\medskip{}

$P_{s}(\overline{u}_{i})=P_{s}^{d}(\overline{\zeta}_{i}),\ \
i=1,2,3.\quad(3.7)$''

\medskip{}

Considering $g:\mathbb{R}\rightarrow\mathbb{R}$ defined by
$g(\varsigma):=\left(2\nu^{-1}
\varsigma+\alpha^{2}\right)(\mu+\varsigma)^{2}$, in fact,
$\overline{\zeta}_{1}(x)$ is the unique solution of the equation
$g(\varsigma)=\beta^{2}(x)$ on the interval $[-\mu,\infty)$, that is
$g(\overline{\zeta}_{1}(x))=\beta^{2}(x)$ and
$\overline{\zeta}_{1}(x)\ge-\mu$, while $\overline{\zeta}_{2}(x)$
and $\overline{\zeta}_{3}(x)$ are the unique solutions of the
equation $g(\varsigma)=\beta^{2}(x)\leq\eta$ on $[\rho,-\mu]$ and
$[-\tfrac{1}{2}\nu\alpha^{2},\rho]$, respectively. We give this
argument later on.

Besides the fact that it is not explained how
$\frac{\sigma(s)+\alpha
\overline{\zeta}_{i}(s)}{\overline{\zeta}_{i}(s)+\mu}$ is defined in
the case $\overline{\zeta}_{i}(s)+\mu=0$ (which is possible if
$\beta(s)=0$) the only mention to $\overline{u}_{i}$ is in the
following paragraph of the proof of \cite[Th.~2]{Gao/Ogden:08}:

\medskip{}

{}``For each solution $\overline{\zeta}_{i}$, $i=1,2,3$, the
corresponding solution $\overline{u}_{i}$ is obtained by rearranging
(2.10) in the form $u_{x}=(\sigma+\alpha\zeta)/(\zeta+\mu)$ and
integrating. For a given $t$, the dead load $\sigma_{1}$ is uniquely
determined by $(3.3)_{3}$. Therefore, there is one
$\overline{u}_{i}(x)$, $i=1,2$ or $3$, satisfying the boundary
condition $\overline{u}_{i}(1)=t$, and this solves (BVP2).''

\bigskip{}

With our reformulation of the problem $(\mathcal{P}_{s})$, in the
statements of \cite[Th.~2, Th.\ 3]{Gao/Ogden:08} one must replace
$\mathcal{U}_{s}$ by $\mathcal{L}^{4}$, $\overline{u}_{i}$ by
$\overline{v}_{i}:=\frac{\sigma+\alpha\overline{\zeta}_{i}}{\overline{\zeta}_{i}+\mu},$
$\overline{u}$ by $\overline{v}$ and $P_{s}$ by $\widehat{P}_{s},$
$\mathcal{U}_{j}$ being a neighborhood of $\overline{v}_{j}$, for
$j=2,3$ (this is possible since the operator
$v\in\mathcal{L}^{4}\rightarrow u= \int_{0}^{x}v\in\mathcal{U}_{s}$
and its inverse $\mathcal{U}_{s}\ni u\rightarrow
v=u_{x}\in\mathcal{L}^{4}$ are linear continuous under the $W^{1,4}$
topology on $\mathcal{U}_{s}$; whence $u\in\mathcal{U}_{s}$ is a
local extrema for $P_{s}$ iff the corresponding
$v\in\mathcal{L}^{4}$ is a local extrema for $\widehat{P}_{s}$).

\medskip{}

We agree that for $\tau^{2}>\eta$ the equation
$\left(2\nu^{-1}\varsigma+\alpha^{2}\right)(\mu+\varsigma)^{2}=\tau^{2}$
has a unique real solution $\varsigma_{1}>-\mu$, while for
$0\leq\tau^{2}\leq\eta$ the preceding equation has three real
solutions $\varsigma_{1},\varsigma_{2},\varsigma_{3}$ with \[
-\tfrac{1}{2}\nu\alpha^{2}\leq\varsigma_{3}\leq\rho\leq\varsigma_{2}\leq-
\mu\leq\varsigma_{1},\]
 where $\rho$ is given in Eq.\ (\ref{rho}).

Indeed, let $g:\mathbb{R}\rightarrow\mathbb{R}$ be defined by
$g(\varsigma):=\left(2\nu^{-1}\varsigma+\alpha^{2}\right)(\mu+\varsigma)^{2}$.
Then $g(\rho)=\eta$ and \[
g^{\prime}(\varsigma)=2\nu^{-1}(\varsigma+\mu)^{2}+2\left(2\nu^{-1}\varsigma+
\alpha^{2}\right)(\mu+\varsigma)=6\nu^{-1}(\varsigma+\mu)(\varsigma-\rho).\]
 The behavior and graph of $g$ are showed in Tables \ref{tab1} and
\ref{tab1-1}.

\begin{table}[h]
 \center \begin{tabular}{c|ccccccccc}
$\varsigma$  & $-\infty$  &  & $-\tfrac{1}{2}\nu\alpha^{2}$  &  &
$\rho$  &  & $-\mu$  &  & $+\infty$\tabularnewline \hline
$g^{\prime}(\varsigma)$  &  & $+$  & $+$  & $+$  & $0$  & $-$  & $0$
& $+$  & \tabularnewline \hline $g(\varsigma)$  & $-\infty$  &
$\nearrow$  & $0$  & $\nearrow$  & $\eta$  & $\searrow$  & $0$  &
$\nearrow$  & $+\infty$\tabularnewline
\end{tabular}

\caption{The behavior of $g$.}

\label{tab1}
\end{table}

\begin{table}[h]
 \center \includegraphics[scale=0.6]{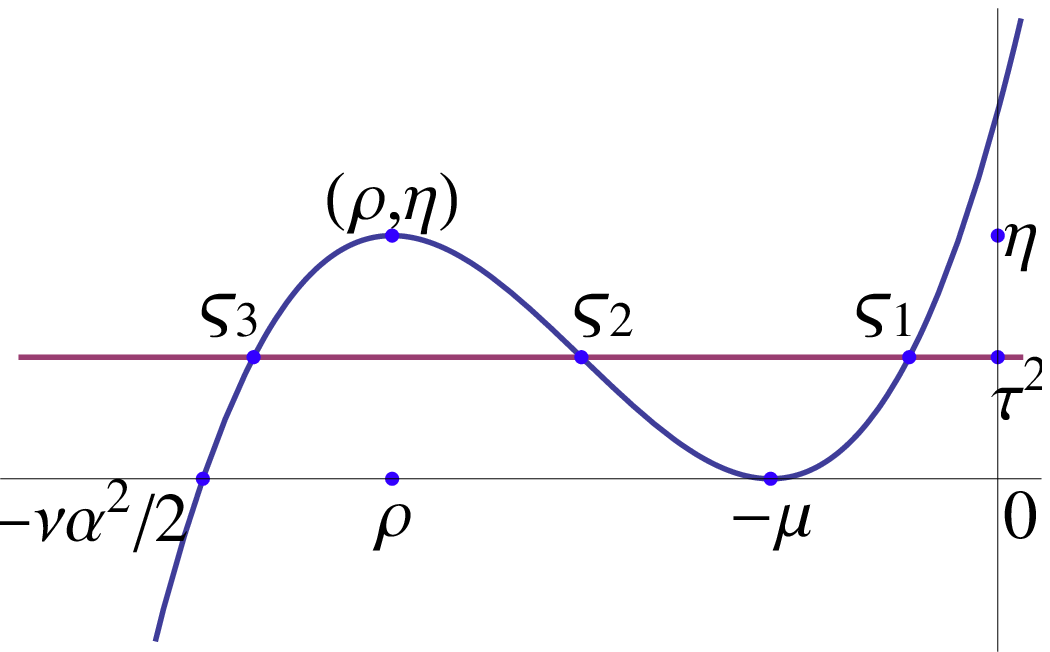}

\caption{The graph of $g$.}

\label{tab1-1}
\end{table}

Note that for $\tau=0$ we have $\varsigma_{1}=\varsigma_{2}=-\mu$,
$\varsigma_{3}=-\tfrac{1}{2}\nu\alpha^{2}$.

For $\tau\in\mathbb{R}$ consider also the function \[
h_{\tau}:\mathbb{R}\setminus\{-\mu\}\rightarrow\mathbb{R},\quad
h_{\tau}
(\varsigma):=-\frac{1}{2}\left[\frac{\tau^{2}}{\varsigma+\mu}+2\alpha\tau
+\alpha^{2}(\varsigma+\mu)+\nu^{-1}\varsigma^{2}\right].\]
 Note that $h_{0}$ is the restriction to $\mathbb{R}\setminus\{-\mu\}$
of the continuous function
$\hat{h}_{0}:\mathbb{R}\rightarrow\mathbb{R}$ defined by
$\hat{h}_{0}(\varsigma):=-\frac{1}{2}\left[\alpha^{2}(\varsigma+\mu)
+\nu^{-1}\varsigma^{2}\right]$; clearly
$\hat{h}_{0}(-\mu)=-\frac{1}{2}\nu^{-1}\mu^{2}$.

Then \[
h_{\tau}^{\prime}(\varsigma)=-\frac{1}{2}\left(-\frac{\tau^{2}}{(\varsigma+\mu)^{2}}
+\alpha^{2}+2\nu^{-1}\varsigma\right)=-\frac{1}{2}\frac{g(\varsigma)-\tau^{2}}
{(\varsigma+\mu)^{2}}\quad\forall\varsigma\in\mathbb{R}\setminus\{-\mu\}.\]
 Taking into account the above discussion (note also the graph of
$g$), the behavior of $h_{\tau}$ is presented in Table \ref{tab2}
for $\tau^{2}>\eta$ and in Table \ref{tab3} for
$0<\tau^{2}\leq\eta$.

\begin{table}[h]
 \center\begin{tabular}{c|ccccccc}
$\varsigma$  & $-\infty$  &  & $-\mu$  &  & $\varsigma_{1}$  &  &
$+\infty$\tabularnewline \hline $h_{\tau}^{\prime}(\varsigma)$  &  &
$+$  & $|$  & $+$  & $0$  & $-$  & $0$\tabularnewline \hline
$h_{\tau}(\varsigma)$  & $-\infty$  & $\nearrow$  &
$^{+\infty}|_{-\infty}$  & $\nearrow$  & $h_{\tau}(\varsigma_{1})$
& $\searrow$  & $-\infty$\tabularnewline
\end{tabular}

\caption{The behavior of $h$ for $\tau^{2}>\eta$.}

\label{tab2}
\end{table}

\begin{table}[h]
 \center\begin{tabular}{c|ccccccccccc}
$\varsigma$  & $-\infty$  &  & $\varsigma_{3}$  &  & $\varsigma_{2}$
& & $-\mu$  &  & $\varsigma_{1}$  &  & $+\infty$\tabularnewline
\hline $h_{\tau}^{\prime}(\varsigma)$  & $+$  & $+$  & $0$  & $-$  &
$0$  & $+$  & $|$ & $+$  & $0$  & $-$  & \tabularnewline \hline
$h_{\tau}(\varsigma)$  & $-\infty$  & $\nearrow$  &
$h_{\tau}(\varsigma_{3})$ & $\searrow$  & $h_{\tau}(\varsigma_{2})$
& $\nearrow$  & $^{+\infty}|_{-\infty}$ & $\nearrow$  &
$h_{\tau}(\varsigma_{1})$  & $\searrow$  & $-\infty$\tabularnewline
\end{tabular}

\caption{The behavior of $h$ for $0<\tau^{2}\le\eta$.}

\label{tab3}
\end{table}


For $\tau=0$ we have that $\hat{h}_{0}$ is increasing on $(-\infty,-
\tfrac{1}{2}\nu\alpha^{2}]$ and decreasing on
$[-\tfrac{1}{2}\nu\alpha^{2},+\infty)$.

\bigskip{}

So, when $\beta^{2}>\eta$ on $(0,1)$ by taking $\tau=\beta(x)$ we
obtain a unique (continuous) solution $\overline{\zeta}$ of
\cite[(2.16)]{Gao/Ogden:08} (with $\overline{\zeta}(x)>-\mu$ for
every $x\in(0,1)$), while for $\beta^{2}\leq\eta$ on $(0,1)$ one
obtains three continuous solutions
$\overline{\zeta}_{1},\overline{\zeta}_{2},\overline{\zeta}_{3}$ of
\cite[(2.16)]{Gao/Ogden:08} satisfying \[
-\tfrac{1}{2}\nu\alpha^{2}\leq\overline{\zeta}_{3}\leq\rho\leq\overline{\zeta}_{2}
\leq-\mu\leq\overline{\zeta}_{1}\text{ on }[0,1].\]

\begin{remark} \label{rem3}In the case $\beta^{2}\leq\eta,$ $\overline{\zeta}_{1},
\overline{\zeta}_{2},\overline{\zeta}_{3}$ are not the only possible
solutions of \cite[(2.16)]{Gao/Ogden:08} with
$\zeta\in\mathcal{L}^{\infty}\subset\mathcal{L}^{2}.$ More
precisely, the general measurable solution
$\zeta:[0,1]\rightarrow\mathbb{R}$ of \cite[(2.16)]{Gao/Ogden:08}
has the form $\zeta(x)=\overline{\zeta}_{j}(x)$ for $x\in B_{j},$
$j=1,2,3$, where $B_{1},B_{2},B_{3}$ are measurable pairwise
disjoint subsets of $[0,1]$ such that $[0,1]=B_{1}\cup B_{2}\cup
B_{3}$.

This shows that none of the $\mathcal{L}^{2}$-solutions of
\cite[(2.16)]{Gao/Ogden:08} is isolated in $\mathcal{L}^{2}$ because
all measurable solutions of \cite[(2.16)]{Gao/Ogden:08} are in
$\mathcal{L}^{\infty}$ and given a measurable solution of
\cite[(2.16)]{Gao/Ogden:08} one can modify it on a sufficiently
small subset (by interchanging the values $\overline{\zeta}_{j}$) so
that it stays still a solution and close enough.\end{remark}

In the sequel we assume that $\beta\neq0$, and so
$\lambda(B_{0}^{c})>0$; the case $\beta=0$ is completely
uninteresting.

\emph{Discussion of \cite[(3.9)]{Gao/Ogden:08}.} Assume that
$\beta^{2}>\eta$ on $(0,1)$. As we have seen above,
$P_{s}^{d}(\zeta)\in\mathbb{R}$ only for $\zeta\in A_{1}\supset
A_{0}$, so considering $\sup_{\zeta\in
\mathcal{S}_{a}}P_{s}^{d}(\zeta)$ in \cite[(3.9)]{Gao/Ogden:08} has
no sense. In the sequel we find sets on which
\cite[(3.9)]{Gao/Ogden:08} holds and then try to further enlarge
them.

In this case the unique solution $\overline{\zeta}$ of
\cite[(2.16)]{Gao/Ogden:08} described above has
$\overline{\zeta}+\mu>0$, and so $E_{\overline{\zeta}}=\emptyset$.
According to Theorem \ref{analysis} (i), (iii)~(b) we have relation
(\ref{equa}) with
$\overline{v}=v_{\overline{\zeta}}=\alpha+\beta/(\mu+\overline{\zeta})$.
This shows that \cite[(3.9)]{Gao/Ogden:08} holds if one replaces
$\max_{\zeta\in\mathcal{S}_{a}}P_{s}^{d}(\zeta)$ by $\max_{\zeta\in
A_{0}, \zeta\ge-\mu}P_{s}^{d}(\zeta)$ (note that $\{\zeta\in
A_{0}\mid\zeta\ge-\mu\}\subset\mathcal{S}_{a}$ because
$\nu\alpha^{2}>2\mu$).

In fact we have that \cite[(3.9)]{Gao/Ogden:08} holds if one
replaces $\max_{\zeta\in\mathcal{S}_{a}}P_{s}^{d}(\zeta)$ by
$\max_{\zeta\in A_{1},\zeta\ge-\mu}P_{s}^{d}(\zeta)$. Indeed,
consider $\zeta\in A_{1}$ with $\zeta\ge-\mu$. Hence $E_{\zeta}$ is
negligible since $B_{0}=\emptyset$; so we may (and do) suppose that
$\zeta$ is finite-valued and $E_{\zeta}=\emptyset$. For $x\in[0,1]$,
from the behavior of $h_{\tau}$ with $\tau=\beta(x)$ (see Table
\ref{tab2}), we obtain that $h_{\beta(x)}(\zeta(x))\leq
h_{\beta(x)}(\overline{\zeta}(x))$, whence \[
P_{s}^{d}(\zeta)=\int_{0}^{1}h_{\beta(x)}(\zeta(x))dx\leq\int_{0}^{1}
h_{\beta(x)}(\overline{\zeta}(x))dx=P_{s}^{d}(\overline{\zeta}).\]

Next we study whether the last equality in
\cite[(3.9)]{Gao/Ogden:08} holds when one replaces
$\max_{\zeta\in\mathcal{S}_{a}}P_{s}^{d}(\zeta)$ by $\max_{\zeta\in
A_{1}^{0}}P_{s}^{d}(\zeta)$, where $A_{1}^{0}:=\{\zeta\in
A_{1}\mid\zeta\geq-\tfrac{1}{2}\nu\alpha^{2}\}$. Unfortunately, that
is not true. Indeed consider $\zeta_{n}(x)=-\mu-\gamma x$ for
$x\in[n^{-1},1]$ and $\zeta_{n}(x)=-\mu-\gamma n^{-1}$ for
$x\in[0,n^{-1})$, where $n\ge1$ and
$0<\gamma<\tfrac{1}{2}\nu\alpha^{2}-\mu$. Clearly
$-\mu-\gamma/n\geq\zeta_{n}\geq-\mu-\gamma>-\tfrac{1}{2}\nu\alpha^{2}$
on $[0,1]$, and so $\zeta_{n}\in A_{1}^{0}$ for every $n\ge1$.
Moreover \[
-\int_{0}^{1}\frac{\beta^{2}}{\zeta_{n}+\mu}dx\geq\int_{1/n}^{1}
\frac{\beta^{2}(x)}{\gamma x}dx\geq\frac{\eta}{\gamma}\ln
n\rightarrow\infty,\]
 which proves that $\sup_{\zeta\in A_{1}^{0}}P_{s}^{d}(\zeta)=+\infty.$

In conclusion \cite[(3.9)]{Gao/Ogden:08} holds if $\mathcal{S}_{a}$
is replaced by anyone of the sets $\left\{ \zeta\in
A_{0}\mid\zeta\ge-\mu\right\} $, $\left\{ \zeta\in
A_{1}\mid\zeta\ge-\mu\right\} $.

Actually the argument above shows that \cite[(3.9)]{Gao/Ogden:08}
holds for $\beta^{2}>0$ on $(0,1)$ if $\mathcal{S}_{a}$ is replaced
by anyone of the sets $\left\{ \zeta\in
A_{0}\mid\zeta\ge-\mu\right\} $, $\left\{ \zeta\in
A_{1}\mid\zeta\ge-\mu\right\} $, $\left\{ \zeta\in
A_{0}\mid\zeta>-\mu\right\} $, $\left\{ \zeta\in
A_{1}\mid\zeta>-\mu\right\} $ with
$\overline{\zeta}=\overline{\zeta}_{1}>-\mu$ (when
$\beta^{2}\le\eta$). The fact that $\overline{v}$ is a minimum point
of $\widehat{P}_{s}$ is confirmed by the fact that
$d^{2}\widehat{P}_{s}
(\overline{v})(h,h)=3\int_{0}^{1}(\overline{\zeta}-\rho)h^{2}dx>0$
for every $h\in\mathcal{L}^{4}\setminus\{0\}$ {[}see
(\ref{d2ps-cr}){]}.

\medskip

\emph{Discussion of \cite[(3.10)]{Gao/Ogden:08}.} Assume that
$\beta^{2}\leq\eta$ on $(0,1)$. As above, if $0<\beta^{2}$ on
$(0,1)$ then \cite[(3.10)]{Gao/Ogden:08} holds if
$\{\zeta\mid\zeta>-\mu\}$ is replaced by anyone of of the sets
$\left\{ \zeta\in A_{0}\mid\zeta>-\mu\right\} $, $\left\{ \zeta \in
A_{1}\mid\zeta>-\mu\right\} $. However, $P_{s}^{d}(\zeta)$ is not
defined for any $\zeta\in\mathcal{L}^{2}$ with $\zeta>-\mu$ so the
previous choices are the only natural ones. Indeed, take
$\zeta(x):=-\mu+x\beta^{2}(x)$ for $x\in(0,1)$; then
$\zeta\in\mathcal{L}^{2}\setminus A_{1}$ and $\zeta>-\mu$ on
$(0,1).$

Assume now that $\lambda(B_{0})>0$ (which happens if $\beta$ is zero
on a nontrivial interval). In this case
$E_{\overline{\zeta}_{1}}=B_{0}$.

Consider $\zeta\in A_{1}$ with $\zeta>-\mu$; hence
$E_{\zeta}=\emptyset\subset B_{0}$. For $x\in B_{0}^{c}$ we have
that $h_{\beta(x)}(\zeta(x))\leq
h_{\beta(x)}(\overline{\zeta}_{1}(x))$ (see Table \ref{tab3}), while
for $x\in B_{0}$, because $h_{0}$ is decreasing on
$[-\tfrac{1}{2}\nu a^{2},+\infty)\setminus\{-\mu\}$, we have that \[
h_{\beta(x)}(\zeta(x))=-\tfrac{1}{2}\left[\alpha^{2}(\zeta(x)+\mu)+\nu^{-1}
\zeta(x)^{2}\right]\leq-\tfrac{1}{2}\nu^{-1}\mu^{2}.\]
 Together with relation (\ref{pdsc2}) applied for $\overline{\zeta}_{1}$,
it follows that $P_{s}^{d}(\zeta)\leq
P_{s}^{d}(\overline{\zeta}_{1})$.

Taking $\varepsilon\in(0,1)$ and
$\zeta_{\varepsilon}(x):=\overline{\zeta}_{1}(x)$ for $x\in
B_{0}^{c}$ and $\zeta_{\varepsilon}(x):=-\mu+\varepsilon$ for $x\in
B_{0}$ we see that $\zeta_{\varepsilon}\in A_{0}\subset A_{1}$
(since $\overline{\zeta}_{1}\in A_{0}$), $\zeta_{\varepsilon}>-\mu$,
and \begin{align*} P_{s}^{d}(\zeta_{\varepsilon}) &
=\int_{0}^{1}h_{\beta(x)}
(\zeta_{\varepsilon}(x))dx=\int_{B_{0}^{c}}h_{\beta(x)}(\overline{\zeta}_{1}
(x))dx+\int_{B_{0}}h_{0}(-\mu+\varepsilon)dx\\
 & =P_{s}^{d}(\overline{\zeta}_{1})-\tfrac{1}{2}[\nu^{-1}\varepsilon^{2}
 +\nu^{-1}(\nu\alpha^{2}-2\mu)\varepsilon]\lambda(B_{0}).\end{align*}
 This implies that $\sup_{\zeta\in A_{1},\zeta>-\mu}P_{s}^{d}(\zeta)=
\sup_{\zeta\in
A_{0},\zeta>-\mu}P_{s}^{d}(\zeta)=P_{s}^{d}(\overline{\zeta}_{1}).$

In the present case {[}that is, $\lambda(B_{0})>0${]}
$\overline{v}_{1}$ is not uniquely determined on $B_{0}$. Taking
$\overline{v}_{1}= v_{\overline{\zeta}_{1}}$, i.e.,
$\overline{v}_{1}(s):=\frac{\sigma+\alpha\overline{\zeta}_{1}}{\overline{\zeta}_{1}+\mu}$
for $s\in B_{0}^{c}$ and $\overline{v}_{1}(s):=\alpha$ for $s\in
B_{0}$ (the natural choice due to the convention $0/0:=0$), we see
from (\ref{pps}) applied for $\zeta=\overline{\zeta}_{1}$ and $v=0$
that $\widehat{P}_{s}(\overline{v}_{1})\neq
P_{s}^{d}(\overline{\zeta}_{1})$, and so \cite[(3.10)]{Gao/Ogden:08}
does not hold.

Again from (\ref{pps}), we see that in order to have that
$\widehat{P}_{s}(\overline{v}_{1})=P_{s}^{d}(\overline{\zeta}_{1})$
we need to have
$\overline{v}_{1}:=v_{\overline{\zeta}_{1}}+\chi_{B_{0}}v$ with
$v\in\mathcal{L}^{4}$ and $v^{2}-\alpha^{2}+2\nu^{-1}\mu=0$ a.e.\ in
$E_{\overline{\zeta}_{1}}=B_{0}$ . In this case, according to
Theorem \ref{analysis} (iii)(a),
$(\overline{v}_{1},\overline{\zeta}_{1})$ is a critical point of
$\widehat{\Xi}$, and so \cite[(3.10)]{Gao/Ogden:08} holds using
(\ref{equa}) if we replace $\max_{\zeta>-\mu}P_{s}^{d}(\zeta)$ by
$\sup_{\zeta\in A_{0},\zeta>-\mu}P_{s}^{d}(\zeta)$ or
$\sup_{\zeta\in A_{1}, \zeta>-\mu}P_{s}^{d}(\zeta).$

Again, in this case
$d^{2}\widehat{P}_{s}(\overline{v}_{1})(h,h)=3\int_{0}^{1}
(\overline{\zeta}_{1}-\rho)h^{2}dx>0$ for every
$h\in\mathcal{L}^{4}\setminus\{0\}$ {[}see (\ref{d2ps-cr}){]} as a
confirmation of
$\widehat{P}_{s}(\overline{v}_{1})=\min_{v\in\mathcal{L}^{4}}\widehat{P}_{s}(v)$.

\medskip

\emph{Discussion of \cite[(3.11)]{Gao/Ogden:08}.} Assume that
$\beta^{2}\leq\eta$ on $(0,1)$. It is easy to show that
$\{\zeta\in\mathcal{L}^{2}\mid\rho <\zeta<-\mu\}\not\subset A_{1}$,
which proves that
$\{\zeta\in\mathcal{L}^{2}\mid\overline{\zeta}_{3}<\zeta<-\mu\}\not\subset
A_{1}$; take for example $\beta^{2}>0$ and
$\zeta(x)=-\mu+\frac{\rho+\mu}{\eta}x\beta^{2}(x)$, $x\in(0,1)$.

This shows that
$\min_{\overline{\zeta}_{3}<\zeta<-\mu}P_{s}^{d}(\zeta)$ in
\cite[(3.11)]{Gao/Ogden:08} does not make sense. Therefore in
\cite[(3.11)]{Gao/Ogden:08} we replace the set
$\{\zeta\in\mathcal{L}^{2} \mid\overline{\zeta}_{3}<\zeta<-\mu\}$ by
$A_{1}^{2}:=\{\zeta\in A_{1}\mid\overline{\zeta}_{3}<\zeta<-\mu\}$.

Here again $B_{0}=E_{\overline{\zeta}_{2}}$. Since
$\overline{\zeta}_{2}(x)$ is the unique minimum point of
$h_{\beta(x)}$ on $[\overline{\zeta}_{3}(x),-\mu)$ for $x\in
B_{0}^{c}$ and $h_{0}$ is decreasing on
$[\overline{\zeta}_{3}(x),-\mu)=[-\tfrac{1}{2}\nu\alpha^{2},-\mu)$
for $x\in B_{0}$, we obtain that for every $\zeta\in A_{1}^{2}$ we
have \[
P_{s}^{d}(\zeta)=\int_{0}^{1}h_{\beta(x)}(\zeta(x))dx=\int_{B_{0}^{c}}
\ldots+\int_{B_{0}}\ldots\ge\int_{B_{0}^{c}}h_{\beta(x)}
(\overline{\zeta}_{2}(x))dx+\int_{B_{0}}\hat{h}_{0}(-\mu)dx=P_{s}^{d}(\overline{\zeta}_{2}).\]
 As above we obtain that $\inf_{\zeta\in A_{1}^{2}}P_{s}^{d}
(\zeta)=P_{s}^{d}(\overline{\zeta}_{2})$ after taking
$0<\varepsilon<-\mu-\rho$ and considering $\zeta_{\varepsilon}\in
A_{1}^{2}$ given by
$\zeta_{\varepsilon}(x):=\overline{\zeta}_{2}(x)$ for $x\in
B_{0}^{c}$ and $\zeta_{\varepsilon}(x):=-\mu-\varepsilon$ for $x\in
B_{0}$.

As seen in the previous discussion (recall also (\ref{pps})), in
order to have
$\widehat{P}_{s}(\overline{v}_{2})=P_{s}^{d}(\overline{\zeta}_{2})$
in \cite[(3.10)]{Gao/Ogden:08} we must take
$\overline{v}_{2}:=v_{\overline{\zeta}_{2}}+\chi_{B_{0}}v$ with
$v\in\mathcal{L}^{4}$ and $v^{2}-\alpha^{2}+2\nu^{-1}\mu=0$ a.e.\ in
$E_{\overline{\zeta}_{2}}=B_{0}.$ With $\overline{v}_{2}$ chosen
this way we have \[
d^{2}\widehat{P}_{s}(\overline{v}_{2})(h,h)=3\int_{0}^{1}(\overline{\zeta}_{2}-
\rho)h^{2}dx\geq0\quad\forall h\in\mathcal{L}^{4}.\]
 However, in general, this $\overline{v}_{2}$ is not a local minimum
point of $\widehat{P}_{s}$. First this is due to the fact that for
$\beta^{2}=\eta$, we have $\overline{\zeta}_{2}=\rho$,
$E_{\overline{\zeta}_{2}} =B_{0}=\emptyset$, and by direct
computation the polynomial that governs $\widehat{P}_{s}$ (i.e.
$\widehat{P}_{s}(v)=\int_{0}^{1}p(v(x))dx$), namely \begin{equation}
p(y):=\tfrac{1}{2}\mu
y^{2}+\tfrac{1}{2}\nu\left(\tfrac{1}{2}y^{2}-\alpha
y\right)^{2}-(\alpha\mu+\beta)y\label{pol-p}\end{equation}
 has $v_{0}:=\alpha+\beta/(\rho+\mu)$ a critical point which is not
a local extremum since $p'(v_{0})=p''(v_{0})=0$,
$p'''(v_{0})=3\nu\beta/(\rho+\mu)\neq0$ and these facts imply that
$v_{0}$ is not a local extremum point for $p$. This implies that
whenever $\beta^{2}=\eta$, $v_{\rho}(x)=v_{0}$, $x\in(0,1)$ is a
critical point but not a local extremum point of $\widehat{P}_{s}$.
Based on the previous facts it is easy to build a counterexample by
taking $\beta$ such that $\beta^{2}=\eta$ on a nonempty open
sub-interval of $[0,1]$. Hence \cite[(3.11)]{Gao/Ogden:08} is not
true even with the correct choice of $\overline{v}_{2}$ and with
$\{\zeta\in\mathcal{L}^{2}\mid\overline{\zeta}_{3}<\zeta<-\mu\}$
replaced by $A_{1}^{2}$ due to the failure of its first equality.

The next natural question is whether
$\overline{v}_{2}=v_{\overline{\zeta}_{2}}+\chi_{B_{0}}v$ with
$v\in\mathcal{L}^{4}$ and $v^{2}-\alpha^{2}+2\nu^{-1}\mu=0$ a.e.\ in
$E_{\overline{\zeta}_{2}}=B_{0}$ is a local minimum point of
$\widehat{P}_{s}$ when $0<\beta^{2}<\eta$ on $(0,1)$ because in this
case $d^{2}\widehat{P}_{s}(\overline{v}_{2})(h,h)>0$ for every
$h\in\mathcal{L}^{4}\setminus\{0\}.$ The answer is still negative as
the next example shows.

\begin{example} Take $\nu:=\mu:=1,$ $\alpha:=3$ and $\beta:=\sqrt{5}$
(a constant function). Note that $\eta=343/27\simeq12.7>\beta^{2}$.
Then the equation $g(\varsigma)=\beta^{2}$ has the solutions
$\varsigma_{1}=(\sqrt{65}-9)/4,$ $\varsigma_{2}=-2$ and
$\varsigma_{3}=-(\sqrt{65}+9)/4.$ Hence $\overline{\zeta}_{2}$ is
the constant function $-2$ and so
$E_{\overline{\zeta}_{2}}=B_{0}=\emptyset$. It follows that
$P_{s}^{d}(\overline{\zeta}_{2})=h_{\sqrt{5}}(-2)=-3\sqrt{5}$ and
$\overline{v}_{2}(x)=v_{\overline{\zeta}_{2}}(x)=y_{0}:=3-\sqrt{5}.$
Moreover, \[
p\left(y_{0}+h\right)-p(y_{0})=\tfrac{1}{8}h^{2}\big(h-2\sqrt{5}+2\big)
\big(h-2\sqrt{5}-2\big)\quad(h\in\mathbb{R}),\]
 where $p$ is the polynomial in (\ref{pol-p}). Consider $\varepsilon\in(0,1)$
and $v:[0,1]\rightarrow\mathbb{R}$ defined by
$v(x):=y_{0}+2\sqrt{5}$ for $x\in[0,\varepsilon]$ and $v(x):=y_{0}$
for $x\in(\varepsilon,1].$ Then $\left\Vert
v-\overline{v}_{2}\right\Vert
_{\mathcal{L}^{4}}=2\sqrt{5}\varepsilon^{1/4}$ and
$\widehat{P}_{s}(v)-\widehat{P}_{s}(\overline{v}_{2})=-10\varepsilon<0,$
which proves that $\overline{v}_{2}$ is not a local minimum of
$\widehat{P}_{s}.$
\end{example}


\emph{Discussion of \cite[(3.12)]{Gao/Ogden:08}.} Assume that
$\beta^{2}\leq\eta$ on $(0,1)$. First note that
$A_{1}^{3}:=\{\zeta\in\mathcal{L}^{2}\mid-\tfrac{1}{2}\nu\alpha^{2}<
\zeta<\overline{\zeta}_{2}\}\subset A_{0}\subset A_{1}$ since
$\overline{\zeta}_{2}\in A_{0}$, and so $P_{s}^{d}(\zeta)$ makes
sense on $A_{1}^{3}$. More precisely, for $\zeta\in A_{1}^{3}$ we
have that \[
\left(\frac{\beta(x)}{\zeta(x)+\mu}\right)^{2}<\left(\frac{\beta(x)}{\overline{\zeta}_{2}(x)
+\mu}\right)^{2}=2\nu^{-1}\overline{\zeta}_{2}(x)+\alpha^{2}\quad\forall
x\in B_{0}^{c}\]
 and $\frac{\beta(x)}{\zeta(x)+\mu}=0$ for $x\in B_{0}$; so $\frac{\beta}{\zeta
 +\mu}\in\mathcal{L}^{4}$,
whence $\zeta\in A_{0}$.

Since $\overline{\zeta}_{3}(x)$ is the maximum point of
$h_{\beta(x)}$ on
$[-\tfrac{1}{2}\nu\alpha^{2},\overline{\zeta}_{2}(x)]$ for $x\in
B_{0}^{c}$ and $h_{0}$ is decreasing on
$[-\tfrac{1}{2}\nu\alpha^{2},-\mu)$ and
$\overline{\zeta}_{3}(x)=-\tfrac{1}{2}\nu\alpha^{2}$ for $x\in
B_{0}$, we obtain similarly that $P_{s}^{d}(\zeta)\leq
P_{s}^{d}(\overline{\zeta}_{3})$ for every $\zeta\in A_{1}^{3}$ or
equivalently $\sup_{\zeta\in A_{1}^{3}}P_{s}^{d}(\zeta)\le
P_{s}^{d}(\overline{\zeta}_{3})$. In a similar manner one can prove
$\sup_{\zeta\in A_{1}^{3}}P_{s}^{d}(\zeta)
=P_{s}^{d}(\overline{\zeta}_{3})$ (see previous discussions).

Since $\overline{\zeta}_{3}$ is not in $A_{1}^{3}$ for those $\beta$
with $\beta^{2}(x)=0$ or $\beta^{2}(x)=\eta$ at some $x\in(0,1)$,
one must replace $\max_{\zeta\in A_{1}^{3}}P_{s}^{d}(\zeta)$ by
$\sup_{\zeta\in A_{1}^{3}} P_{s}^{d}(\zeta)$. This time
$\widehat{P}_{s}(\overline{v}_{3})=P_{s}^{d}(\overline{\zeta}_{3})$
because $E_{\overline{\zeta}_{3}}=\emptyset$. However, as previously
seen for \cite[(3.11)]{Gao/Ogden:08}, in general $\overline{v}_{3}$
is not a local maximum point of $\widehat{P}_{s}$. So
\cite[(3.12)]{Gao/Ogden:08} is not true under the hypotheses of
\cite[Th.~3]{Gao/Ogden:08} again because its first equality does not
hold.


\section{Conclusions}
\begin{itemize}
\item The statement of \cite[Th.~3]{Gao/Ogden:08} is ambiguous because
$P_{s}^{d}(\zeta)$ is not defined for all $\zeta$ to which it is
referred and $\overline{u}_{1}$ and $\overline{u}_{2}$ are not
clearly and properly defined.
\item The left equalities in \cite[(3.11)]{Gao/Ogden:08} and \cite[(3.12)]{Gao/Ogden:08}
are not true in general even when proper choices are considered for
the sets where the maximization or minimization of $P_{s}$ happens
and correct choices of $\overline{u}_{i}$ are taken.
\item For proper choices of the sets where the maximization or minimization
of $P_{s}^{d}$ is considered, the right equalities in relations
(3.9)--(3.12) of \cite[Th.~3]{Gao/Ogden:08} follow by very
elementary arguments.
\item Note that in Gao's book \cite[page 140]{Gaoo-book} it is said: {}``For
any given critical point
$(\overline{u},\overline{\varsigma})\in\mathcal{L}_{c}$, we let
$\mathcal{U}_{r}\times\mathcal{T}_{r}$ be its neighborhood such
that, on $\mathcal{U}_{r}\times\mathcal{T}_{r},$
$(\overline{u},\overline{\varsigma})$ is the only critical point of
$L$. The following result is of fundamental
importance in nonconvex analysis.\textbf{}\\
\textbf{ Theorem 3.5.2 (Triality Theorem)} Suppose that
$(\overline{u}, \overline{\varsigma})\in\mathcal{L}_{c}$, and
$\mathcal{U}_{r}\times\mathcal{T}_{r}$ is a neighborhood of
$(\overline{u},
\overline{\varsigma})...$''\textbf{}\\
 We think that such a result was used for proving \cite[Th.~3]{Gao/Ogden:08}.
Taking into account Remark \ref{rem3}, we see that, for
$\beta^{2}\leq\eta$, $\widehat{\Xi}$ has no isolated critical
points; hence the previous theorem cannot be used as an argument for
\cite[Th.~3]{Gao/Ogden:08}. Having in view this situation, it would
be interesting to know the precise result the authors used to derive
\cite[Th.~3]{Gao/Ogden:08}. \end{itemize}

\medskip
{\bf Acknowledgement.} The paper was submitted to ``The Quarterly
Journal of Mechanics and Applied Mathematics'' in February 2010
under the title `On a result about global minimizers and local
extrema in phase transition'. Besides the title, the only difference
is that in the Introduction instead of

$(\mathcal{P}_{s})$ : $\min\limits
_{u\in\mathcal{U}_{s}}{\displaystyle  \left\{
P_{s}(u)=\int_{0}^{1}\Bigl[\tfrac{1}{2}\mu
u_{x}^{2}+\tfrac{1}{2}\nu\left(\tfrac{1}{2}u_{x}^{2}-\alpha
u_{x}\right)^{2}\Bigr]dx-F(u)\right\} }$, $\quad$(3.2)''

\noindent there was

$(\mathcal{P}_{s})$ : $\min\limits
_{u\in\mathcal{U}_{s}}{\displaystyle  \left\{
P_{s}(u)=\int_{0}^{1}\Bigl[\tfrac{1}{2}\mu
u_{x}^{2}+\tfrac{1}{2}\nu\left(\tfrac{1}{2}u_{x}^{2}-\alpha
u_{x}\right)^{2}dx-F(u)\Bigr]\right\} }$, $\quad$(3.2)''.


\begin{thebibliography}{7}
\bibitem{Gao:98}D. Y. Gao, Duality, triality and complementary extremum
principles in non-convex parametric variational problems with
applications, \textit{IMA J. Appl. Math.} \textbf{61} (1998)
199--235.

\bibitem{Gao:99}D. Y. Gao, General analytic solutions and complementary
variational principles for large deformation nonsmooth mechanics,
 \textit{Meccanica} \textbf{34} (1999) 169--198.

\bibitem{Gao:00}D. Y. Gao, Analytic solutions and triality theory
for non-convex and nonsmooth variational problems with applications,
 \textit{Nonlinear Anal.} \textbf{42} (2000) 1161--1193.

\bibitem{Gaoo-book}D. Y. Gao,  \textit{Duality Principles in Nonconvex Systems:
Theory, Methods and Applications} (Kluwer, Dordrecht 2000).

\bibitem{Gao/Ogden:08}D. Y. Gao, R. W. Ogden, Multiple solutions
to non-convex variational problems with implications for phase
transitions and numerical computation,  \textit{Quart. J. Mech.
Appl. Math.} \textbf{61} (2008)  497--522.

\bibitem{Gao/Strang:89}D. Y. Gao, G. Strang, Geometric nonlinearity:
Potential energy, complementary energy, and the gap function,
 \textit{Quart. Appl. Math.} \textbf{47} (1989) 487--504.

\bibitem{roy-88}H. L. Royden,  \textit{Real analysis (3rd edition)} (Macmillan
Publishing Company, New York 1988).
\end{thebibliography}
\end{document}